\documentclass[12pt]{article}
\usepackage{amsfonts, amssymb, amsmath, latexsym}

\newcommand{\la}{\langle}
\newcommand{\ra}{\rangle}

\newcommand{\goth}{\mathfrak}

\newcommand{\pz}{\partial}
\newcommand{\C}{\mathbb{C}}
\newcommand{\R}{\mathbb{R}}
\newcommand{\bz}{\bar{\partial}}

\newcommand{\Pz}{\frac{\partial}{\partial z}}
\newcommand{\Bz}{\frac{\partial}{\partial \bar{z}}}
\newcommand{\n}{{\bf n}}
\newcommand{\bb}{\mathbb} 
\newtheorem{theorem}{Theorem}[section]
\newtheorem{proposition}{Proposition}[section]
\newtheorem{corollary}[theorem]{Corollary}
\newtheorem{lemma}[theorem]{Lemma}
\newtheorem{remark}{Remark}[section]

\begin{document}

\title{ Spacelike Surfaces  in  De Sitter $3$-space \\ and their twistor lifts } 
\author{Eduardo Hulett
 \footnote{Partially supported by ANPCyT, CONICET and SECYT-UNC,
Argentina.}\\ email: hulett@famaf.unc.edu.ar  \\
phone/fax:  +54 351 4334052/51\\
C.I.E.M. - Fa.M.A.F.\\
 Universidad Nacional de C\'ordoba,\\
Ciudad Universitaria,\\
5000 C\'ordoba, Argentina. }
\date{}
\maketitle

 \begin{abstract}
 We deal here with the geometry of the so-called twistor fibration $\mathcal{Z} \to \bb{S}^3_1$ over  the De Sitter $3$-space, where the total space   $\mathcal{Z}$  is a five dimensional reductive homogeneous space with two canonical invariant almost CR structures. Fixed the  normal metric on $\mathcal{Z}$ we study   
  the harmonic map equation for smooth maps of Riemann surfaces into $\mathcal{Z}$. A characterization  of spacelike surfaces with harmonic   twistor lifts to $\mathcal{Z}$ is obtained. Also it is shown that the harmonic map equation for twistor lifts can be formulated as the curvature vanishing  of an $\bb{S}^1$-loop of connections i.e.  harmonic twistor lifts exist within $\bb{S}^1$-families.  Special harmonic maps such as holomorphic twistor lifts  are also considered and some remarks concerning (compact) vacua of the twistor energy are given.  
 \end{abstract}
Math. Subject Classification (2000) 53C43, 53C42, 53C50.

\section{Introduction} 
J. Eells and S. Salamon~\cite{eells-salamon}  were able to obtain  conformal harmonic maps of Riemann surfaces into a $3$-dimensional Riemannian manifold $N$ by projecting Cauchy-Riemann holomorphic maps with values in the unit tangent bundle $T^1 N$. Replacing the riemannian $3$-manifold $N$ by a Lorentzian $3$-manifold $L$, it is then natural to look for  analogous results.
The "unit" tangent bundle $T^{1}L =\{ x \in TL : \| x \|^2 = \pm 1 \}$, of a Lorentzian $3$-manifold $L$ splits into the disjoint union of the timelike unit tangent bundle $T_{-1}^1 L=\{ x \in TL : \| x \|^2 = -1 \} $ and the spacelike unit tangent bundle $T_{+1}^1 L=\{ x \in TL : \| x \|^2 = +1 \}$. If $L$ comes equipped with  a global notion of time orientedness, then    
any connected spacelike surface in $L$ has a uniquely defined future oriented timelike unit normal field $\n$  which may be viewed as a kind of Gauss lift of the surface to the subbundle $\mathcal{Z} \subset T_{-1}^1 L$  consisting of all future oriented timelike  unit tangent vectors of $L$.  We call $\mathcal{Z}$ the twistor bundle of the Lorentzian $3$-manifold $L$, by analogy with the construction of  hyperbolic twistor spaces due to D. Blair, J. Davidov and O. Muskarov~\cite{blair-davidov-muskarov}.\\

    Our goal  in this paper is to give a detailed account of the  geometry of the twistor bundle $\mathcal{Z}$ constructed over De Sitter $3$-space or pseudosphere $\bb{S}^3_1$. Also we intend to understand the geometry behind  the harmonic map equation for maps from Riemann surfaces into $\mathcal{Z}$, when the normal metric is fixed on the target. In particular we characterize conformally immersed (hence spacelike) surfaces in   $\bb{S}^3_1$ whose lifts are harmonic maps.     \\

     As a homogeneous manifold of the simple Lie group $SO_o(3,1)$  the twistor space $\mathcal{Z}$ is a (non-compact) reductive quotient  equipped with an invariant horizontal distribution $\goth{h} \subset T \mathcal{Z} $ which is  just  the orthogonal complement with respect to the normal metric to the vertical distribution  of the principal bundle $\mathcal{Z} \to G^+_2(\R^4_1)$, where the last is the Grassmann manifold of oriented spacelike subspaces of $\R^4_1$.  The distribution $\goth{h} \subset T \mathcal{Z}$ 
     supports two invariant almost complex structures  which arise from the natural complex structure on the fibers of $\mathcal{Z} \to \bb{S}^3_1$ (diffeomorphic to the hyperbolic $2$-space $\bb{H}^2$), and the almost complex structure naturally attached to every  point of  $\mathcal{Z}$. These almost complex structures give rise to corresponding CR structures on $\mathcal{Z}$.  
Since the horizontal distribution  $\goth{h}$ is contact, its integral manifolds have dimension at most two. In fact,   twistor lifts of conformally immersed Riemann surfaces provide  examples of Legendrian manifolds (maximal integral  manifolds) of the distribution $\goth{h}$.  
   We consider the energy of smooth maps $\phi:M \to \mathcal{Z}$ determined by the normal metric on  $\mathcal{Z}$. In this context we  investigate those lifts which are critical points of the energy i.e. which satisfy the harmonic map equation.  
   Our main result is Theorem~\ref{harmonictwsitorlift}   which gives a characterization conformally immersed surfaces in $\bb{S}^3_1$ with harmonic twistor lift. Also in Theorem~\ref{S^1deformations} we show that the harmonic map equation for  twistor lifts is a completely integrable system i.e. harmonic twistor lifts exist within $\bb{S}^1$-families.  \\

The paper is organized as follows. In the first section we derive the basic  structure equations of conformally immersed  Riemann  surfaces  in $\bb{S}^3_1$. 
The second section is devoted to  study  the geometry of the twistor fibration  $\mathcal{Z} \to \bb{S}^3_1$.  We introduce the  horizontal distribution  $ \goth{h} \subset T \mathcal{Z} $ and define two invariant almost complex  structures $J', J''$  on $\mathcal{Z}$. 
In the third  section we  derive the harmonic map equation for smooth maps of Riemann surfaces with values in $\mathcal{Z}$ in terms   of the Maurer Cartan one form $\beta$ of the reductive space $\mathcal{Z}$.  
In section 4 we  characterize  spacelike surfaces with  harmonic twistor lifts.   
In section 5 we deal with one (complex) parameter deformations of harmonic twistor lifts. Although  $\mathcal{Z}$ is not a symmetric space, we show that the harmonic map equation for twistor lifts can be formulated as a  loop of flat connections, a characteristic property of integrable systems (see~\cite{burstall-pedit} for instance).  We also consider special twistor lifts such as holomorphic ones. 
However, we have not dealt here with specific calculations using loop groups techniques to produce examples of harmonic lifts. This will be considered in  another paper. Finally  in the last section we compute  the energy of twistor lifts and establish a relationship with the Willmore energy. Some remarks are given concerning compact vacua of the twistor energy.

\section{Spacelike surfaces in $\bb{S}^3_1$ } 
 Denote by $\R^4_1$ the   real $4$-space $\R^4$ with coordinates $(x_1, x_2, x_3, x_4)$  equipped with the Lorentz metric 
 $$\la .,.\ra=dx^2_1 + dx^2_2 + dx^2_3 - dx^2_4.$$\\
 \medskip
 De Sitter $3$-space    is defined as the unit sphere  in $\R^4_1$ 
  $$
      \bb{S}^3_1 = \{ x \in  \R^4_1 : \langle x,x \rangle =1 \},$$
 on which the ambient Lorentz metric induces  a pseudometric   $\langle .,. \rangle $ with signature $(++-)$, and so it  becomes a Lorentz $3$-manifold with  constant curvature one.\\
  Let  $I_{3,1} = diag (1,1,1,-1)$,  then the  simple Lie group $SO_o(3,1)= \{A \in Gl_4 (\R) : A^t I_{3,1} A = I_{3,1}, a_{44}>0 \}$ acts transitively on  $\bb{S}^3_1$ by isometries. 
  A global time orientation of $\bb{S}^3_1$ is obtained by declaring  a timelike vector $X \in T_x \bb{S}^3_1$ to be  {\it future-pointing} if $ \la X, V_x \ra < 0$, where $V$ is the unit timelike Killing vector field    $V$ on $\bb{S}^3_1$ given by 
$$   V_x = \frac{d}{dt} \left|_{t=0} \exp(t X_0).x,\right. \,\,\, x \in \bb{S}^3_1,$$
and 
 \begin{equation*}
  X_0 = \begin{pmatrix}    0 & 0 & 0 & 1 \\
                     0 & 0 & 0 & 0 \\
                     0 & 0 & 0 & 0\\
                     1 & 0 & 0 & 0\\
                     \end{pmatrix} \in \mathfrak{so}(3,1).
   \end{equation*} 
It is easily seen that  a timelike vector  $X \in T_x \bb{S}^3_1$ is future pointing if and only if  after parallel translation to the origin of $\R^4_1$ it satisfies $X_4 >0$.\\ 
 
 An immersion $f : M \to \bb{S}^3_1$  of a Riemann surface is   conformal if  
   $\la 
  f_z ,  f_z \ra^{c}=0$, for every  local complex
  coordinate $z=x+iy$ on $M$, where 
  $$\Pz = \frac{1}{2} ( \frac{\partial}{\partial x}- i \frac{\partial}{\partial y}), \quad \Bz = \frac{1}{2} ( \frac{\partial}{\partial x} + i \frac{\partial}{\partial y}),$$
are the  complex partial derivatives  and  $\la \, , \, \ra^{c}$ is the complex bilinear  extension of
   the Lorentz metric to $\C^4$:
   $$  \la z,w \ra^c =  z_1 w_1 + z_2 w_2 + z_3 w_3 - z_4 w_4.  $$
   
    Thus $ f$   is conformal if and only if    on any   local complex coordinate $z =x + iy$ it satisfies 
   \begin{equation} \label{confID}
       \la f_x , f_y \ra =0, \quad \| f_x \|^2 =\| f_y \|^2.
  \end{equation}
 
In particular \eqref{confID} implies that  every conformal immersion $f:M \to \bb{S}^3_1$
  is space-like  i.e.  the induced metric $g=f^* \la \,,\,
  \ra$ is positive definite or riemannian on $M$. \\
  
 Let $f:M \to \bb{S}^3_1$ be a conformal immersion of a connected Riemann surface $M$ and  $\n : M \to T \bb{S}^3_1$  a smooth  future-pointing unit timelike vector field along $f$ which is normal to the immersed surface at each point. Such vector field exists by  the time orientation of   $\bb{S}^3_1$.         
 We fix on $M$ the induced Riemannian metric $g= f^* \la \,,\, \ra$ so that  $ f: (M,g) \to \bb{S}^3_1$ is a spacelike isometric immersion. 
    The {\it $2$nd Fundamental form} of the space-like surface $f:M \to \bb{S}^3_1$  is given by 
    $$ II  = - \la df, d \n \ra. $$
    On any local  chart $(U, z=x+iy)$ of $M$ we introduce   
          a conformal parameter $u$ defined by $\la \pz f, \pz f\ra=e^{2u}$, so that    $g|_U=2 e^{2u}( dx^2 +   dy^2)$. 
 The {\it Mean curvature} of $f$ is  defined by $H= \frac{1}{2}trace II $, which   in terms of $f$ and $u$ is   given by 
      $$ H= -e^{-2u} \la  f_{\bar{z} z}, \n \ra. $$
    Since $f$ is conformal we have
    \begin{equation*}
    \begin{array}{l}
    2\la  f_{\bar{z} z} ,  f_z \ra^{c} = \Pz \la  f_z,  f_z
    \ra^{c} =0\\
    2\la  f_{\bar{z} z} , f_{\bar{z} } \ra^{c} = \Pz \la f_{\bar{z} }, f_{\bar{z} }
    \ra^{c} =0,\\
    \end{array}
    \end{equation*}
    hence     $f_{\bar{z} z}$ has no tangential component. We obtain the 
      structural equations of $f$:  
     \begin{equation*}
     \begin{array}{ll} 
         f_{\bar{z} z} = -e^{2u} f + e^{2u}H. \n. & (i)\\
        f_{z z} = 2  u_z .  f_z + \xi. \,\n, & (ii)\\
        \,  \n_z = H.  f_z + e^{-2u} \xi.  f_{\bar{z}}, & (iii)\\
       \end{array}
       \end{equation*}
       where  $q := \xi \, dz \otimes dz = - \la
 f_{z z}, \n \ra^{c} dz^2$ is the Hopf quadratic complex differential. Zeros  of  $q$ are the  umbilic points of $M$. 
We say that  $f$ is isotropic or totally umbilic if and only if  $\xi \equiv 0$. 
If $H \equiv 0$ then the conformal immersion $f$ is harmonic or maximal. In this case $f$ satisfies 
\begin{equation} \label{harmapeq}
  f_{\bar{z} z} = -\la f_z ,  f_{\bar{z} } \ra^{c}f,  
  \end{equation} 
  where   $\la  f_z,  f_{\bar{z} } \ra^{c} = e^{2u}$.  
Away from  umbilic points of $f$, from equation (ii) we obtain    
\begin{equation} \label{normalvector}
\n = \frac{1}{\xi}. ( 2  u_z f_z -  f_{zz}),
\end{equation}
which allows us to  recover the normal vector field from the immersion.\\
The compatibility conditions of the  structure equations  are
 the Gauss-Codazzi equations:
\begin{equation} \label{gauss-codazzi}
 \begin{array}{ll}
  2  u_{\bar{z} z} = (H^2-1)e^{2u} - |\xi|^2 e^{-2u} & \text{(Gauss)}, \\
  \xi_{\bar{z} } = e^{2u}  H_z & \text{(Codazzi)}.  \\
\end{array}
 \end{equation} 
Conversely, it is known that any solution of these equations defines a surface in $\bb{S}^3_1$ up to an isometry.
 From Codazzi's eq. a  surface has contant mean curvature  $H$  if and only if   $\xi$ is 
holomorphic. Hence   $H= const.$ and  $\xi \not \equiv 0$, then the  umbilic points are isolated. \\

For a conformal immersion $f:M \to \bb{S}^3_1$, we consider the induced metric 
$g=f^* \la \,,\, \ra$ on $M$, hence $f:(M,g) \to \bb{S}^3_1$ is an isometric space-like immersion.  In terms of the conformal parameter $u$,   the induced metric is given by $g=2 e^{2u}dz \otimes d\bar{z}$, and 
the Gaussian curvature of the surface $(M,g)$ is just the curvature of the induced metric and  is given by 
$$K= - \Delta_g u = -2 e^{-2u}  u_{\bar{z} z}, $$
 where $\Delta_g$ is the Laplace operator on $M$ determined by $g$.  In complex coordinates $\Delta_g=  2 e^{-2u} \Bz \Pz$. 
Thus Gauss equation becomes   
\begin{equation}\label{gausseq}
K=1-H^2 + |\xi|^2 e^{-4u}.
\end{equation}
Here $\| q \| = |\xi|^2 e^{-4u}$ is the intrinsic norm of the Hopf differential. If $\lambda_1, \lambda_2$ are the principal curvatures of the immersed surface, then it follows that   
\begin{equation}\label{principalcurvatures}
|\xi|^2 e^{-4u} = \frac{1}{4}(\lambda_1 - \lambda_2)^2,
\end{equation}
so that Gauss equation reads
\begin{equation}\label{gausseq2}
    K+H^2-1 = \frac{1}{4}(\lambda_1 - \lambda_2)^2.
    \end{equation}

\section{The twistor bundle  of $\bb{S}^3_1$ and its geometry}

 By definition the fiber $G_v$ of the Gauss bundle $G \to  \bb{S}^3_1$    over a  point  $v \in \bb{S}^3_1$ is the totality  of  $2$-planes passing through the origin  of $T_v \bb{S}^3_1$. 
 If one is interested only in the geometry  of  spacelike surfaces in $ \bb{S}^3_1$, the Gauss bundle is too big. Restricting  the fibers  of $G$  by allowing only spacelike subspaces, one obtains a  new bundle. We define the twistor bundle $\mathcal{Z} \subset G$ of $\bb{S}^3_1$ by demanding that the fiber $\mathcal{Z}_v$ over a point $v \in \bb{S}^3_1$ be the set  of spacelike $2$-planes in  $T_v \bb{S}^3_1$.  A spacelike $2$-plane $V \subset T_v \bb{S}^3_1$ has  two timelike unit normal vectors  of which only  one  is (according to our definition) future pointing and we choose  it to  fix  the desired orientation on $V$. Hence any $V \in \mathcal{Z}_v$ determines and  is determined by a unit timelike future pointing vector $w \in T_v \bb{S}^3_1$  by requiring $w^{\bot} = V \subset T_v \bb{S}^3_1$. Translating $w$ and $V$    to the origin of $\R^4_1$, they satisfy $\la w,w \ra =-1$, $\la v, w \ra =0$ with $w_4 >0$ and $V =  [v \wedge w ]^{\bot}$. Note that  
  $w$ defines a point in the upper half $\bb{H}^3_+$  of the hyperboloid $\{ x \in \R^4_1: \la x,x \ra =-1 \}$, which is the unbounded realization of the three dimensional real hyperbolic space. Hence the total space of the twistor bundle of $\bb{S}^3_1$ is just  
  \begin{equation}\label{thetwistorbundle}
  \mathcal{Z} = \{ (v,w) \in \bb{S}^3_1 \times \bb{H}^3_+ \subset \R^4_1 \times \R^4_1 :   \la v,w \ra =0 \},
  \end{equation}
   where the projection map $\pi : \mathcal{Z} \to \bb{S}^3_1$ is simply  $\pi (v,w)=v$. 
The  fiber  $\mathcal{Z}_v = \pi^{-1}(v)$ over  $v \in \bb{S}^3_1$ identifies with  
$    v^{\bot} \cap \bb{H}^3_+ ,$
 which is a copy of hyperbolic $2$-space $\bb{H}^2$ totally geodesic immersed  
in $\bb{H}^3_+$, hence a complex manifold.\\
Note that a second fibration $\pi'' : \mathcal{Z} \to \bb{H}^3_+$ is obtained by projection on the second factor, $\pi'': (v,w) \mapsto w$. It is clear that
the fiber of an element $ w \in  \bb{H}^3_+$ is the $2$-sphere $ w^{\bot} \cap \bb{S}^3_1$.

\begin{remark} Any   $w \in \mathcal{Z}_v$ determines the   oriented  spacelike $2$-plane $V=[v \wedge w]^{\bot} \subset T_v \bb{S}^3_1$  to which there is associated   
    the  rotation   $J_w :V \to V$ of angle $\frac{\pi}{2}$ compatible with the orientation on $V$. So that  identifying  
$$w \leftrightarrow [v \wedge w]^{\bot} \leftrightarrow J_w , $$ 
we may think of $\mathcal{Z}_v \equiv \bb{H}^2 $ as the set of all oriented spacelike planes in $T_v \bb{S}^3_1$ and also, as the set of all oriented rotations in $T_v \bb{S}^3_1$. One may view a section of $\mathcal{Z}$ as a field of rotations $ \{ v \mapsto J_v \}$, or  equivalently as a distribution $\mathcal{D}$ of oriented spacelike $2$-planes in $\bb{S}^3_1$.  
\end{remark}

Let $M$ be a connected Riemann surface and $f:M \to \bb{S}^3_1$ a conformal immersion. It is  not hard to  show that there exist a uniquely defined  future-oriented unit normal vector field $\widehat{f}$   along $f$ satisfying  
$$    \widehat{f}(x) \in T_{f(x)} \bb{S}^3_1, \quad \widehat{f}(x) \bot df_x(T_x M), \forall x \in M.$$   
Thus the field $\widehat{f}$ determines a  unique smooth map $\n :M \to \bb{H}^3_+ $ such that
\begin{equation}\label{definitionTwistorLift}
              \widehat{f}(x) = (f(x), \n(x)) \in \mathcal{Z},\,\, x \in M.
\end{equation}
We call  $\widehat{f} :M \to \mathcal{Z}$ the {\it twistor lift} of the conformal immersion $f$, and $\n : M \to \bb{H}^3_1$ its  {\it normal Gauss map}.\\

Since $\mathcal{Z}$ is a submanifold of $ \R^4_1 \times \R^4_1$,  the tangent space of $\mathcal{Z}$ at $(v,w)$  is  given by  
$$
 T_{(v,w)} \mathcal{Z} = \{ (x,y) \in  \R^4_1 \times \R^4_1: \la x,v \ra =
  \la y,w \ra =0,\,\, \la x,w \ra + \la v,y \ra =0 \}.
$$
Hence fixed the base point $o=(e_1, e_4) \in \mathcal{Z}$ then   
  $(x,y) \in T_o \mathcal{Z}$ if and only if 
$$x=( 0 , x_2,
x_3, x_4)^T ,\quad y=(x_4, y_2, y_3  , 0)^T.$$
 Then   we may   identify
\begin{equation}\label{identif1}  \goth{p}:=T_o \mathcal{Z} \ni (x,y) \equiv \begin{pmatrix}
 0     &   x_2  &   x_3  &   x_4 \\
 -x_2  &   0    &   0    &   y_2\\
 -x_3  &   0    &   0    &   y_3\\
 x_4   &   y_2  &   y_3  &   0   \end{pmatrix} \in \goth{so}(3,1). 
 \end{equation}
 
 On the other hand the transitive 
 left action of  $SO_o(3,1)$ on $\mathcal{Z}$ given by 
$    g.(v,w) = (g.v, g.w) $
allows to identify  $\mathcal{Z}$ with 
   the quotient $ \displaystyle{ {SO_o(3,1)/ K}}$,  where 
\begin{equation} \label{subgroupK}
     K= \{ diag (1, A, 1), A \in SO(2) \},
\end{equation}
is  the isotropy subgroup of the base point $o \in \mathcal{Z}$.\\ 
 Decompose 
  $\goth{so}(3,1)= \goth{k}
    \oplus \goth{p}$,
     where
$$
 \goth{k}= \{        \begin{pmatrix}
        0 & 0 & \,\,0 & 0\\
        0 & 0 & -a & 0 \\
        0 & a & \,\,0& 0 \\
        0 & 0 & \,\,0 & 0\\
        \end{pmatrix}, a \in R \},$$
     is the Lie algebra of
    $K  \simeq SO(2)$. 
Then the decomposition   is reductive since $\goth{p} \equiv T_o \mathcal{Z}$ defined above is  $Ad(K)$-invariant.
From now on  we choose the    $Ad(K)$-invariant  inner product $\la .,. \ra$ on  $\goth{p} \equiv T_o \mathcal{Z}$ defined by     
\begin{equation}\label{inn.prod.so(3,1)}
 \la A,B \ra = -\frac{1}{2}\, trace(A.B), \, \,\, A, B \in \goth{p},
 \end{equation}
Note that     $\|A\|^2=
x^2+y^2-c^2-z^2-w^2$,  $ \forall A \in
\goth{p}$,  
hence~\eqref{inn.prod.so(3,1)}  gives rise to  an
$SO_o(3,1)$-invariant pseudo-metric on $\mathcal{Z}$ of
signature $(++---)$ denoted also by $\la .,. \ra$, the so-called {\it normal metric}.  Since  $\la . , . \ra$ is the restriction of (a multiple of ) the Killing form of $\goth{so}(3,1)$ to $\goth{p} \times \goth{p}$,  $( \mathcal{Z}, \la .,. \ra )$ is naturally reductive. In this case the natural projection $SO_o(3,1) \to \mathcal{Z}$ is a (pseudo) riemannian submersion, in which the bi-invariant (pseudo) metric induced by the Killing form is considered on $SO_o(3,1)$. From now on we shall consider on $\mathcal{Z}$ only the normal metric.  

\begin{remark}  Another remarkable $SO_o(3,1)$-invariant  metric is the one that makes   the inclusion  $\mathcal{Z} \subset \R^4_1 \times \R^4_1$ an isometric immersion.   
\end{remark}

 Given    $(v,w)\in \mathcal{Z}$  the oriented
 space-like $2$-plane  
 $V = [v \wedge w ]^{\bot}$ defines a  point of  the Grassmannian $G_2^+ (\R^4_1)$ of all oriented spacelike planes in $\R^4_1 $. Definining the  projection  map $\pi':  \mathcal{Z} \to G_2^+
(\R^4_1)$,  by 
$$             \pi'(v,w) = [ v \wedge
  w]^{\bot},
$$
 we obtain  an $SO_o(1,1)$-principal bundle  $\pi':  \mathcal{Z} \to G_2^+
(\R^4_1)$ where  the right action of
$SO_o(1,1)$ on the total space $\mathcal{Z}$ is given by

\begin{equation}\label{rightaction}
(v,w).\begin{pmatrix}
 \cosh(t) & \sinh(t) \\
 \sinh(t) & \cosh(t)\\
 \end{pmatrix} = 
 (\cosh (t) v + \sinh(t) w , \sinh (t) v + \cosh(t)
 w).
\end{equation}
Pick a point  $V \in G_2^+(R^4_1)$ and let  $\{v,w \}$ be  an oriented basis of of $V^{\bot}$ with  $\|v \|^2 = -\|w \|^2 =1, \, \la
v, w \ra =0$ and $w_4 >0$, so that $(v,w) \in \mathcal{Z}$ and $V=[v \wedge w]^{\bot}$. Define the curve $\tau(t)=\cosh(t)v+\sinh(t)w \subset V^{\bot}= [v \wedge w]$, then  $\gamma(t)=( \tau(t) , \tau'(t)) \in \mathcal{Z} $ for all $t \in \R$ and $\gamma(0) = (v, w)$.  We  conclude that   
$$  \pi'^{-1}(V)= \{ ( \tau(t) , \dot{\tau}(t)):
t \in \R \}.
$$
 In other words~\eqref{rightaction} shows that  the $\pi'$-fibre 
            through $(v,w) \in \mathcal{Z}$ is just the (right) orbit   
           $$ (v,w).SO_o(1,1)= \{ (v,w). \exp(t Z_0): t \in \R \}, \quad 
Z_0 :=
\begin{pmatrix}
            0&1 \\
           1& 0 \\
           \end{pmatrix} \in \goth{so}(1,1).$$

Note that for any  $(v,w) \in \mathcal{Z}$ we have  
$$\frac{d}{dt}|_{t=0} (v,w).\exp(t Z_0)= \frac{d}{dt}|_{t=0} (\tau(t), \dot{\tau}(t)) = (w,v) \in T_{(v,w)} \mathcal{Z}. $$

This suggests defining the characteristic (or Hopf) vector field $h$ on $\mathcal{Z}$  by 
\begin{equation}\label{hopfvectorfield}
 h_{(v,w)} = (w,v) \in T_{(v,w)}\mathcal{Z}.
 \end{equation}
Hence  $h$ is a unit timelike vector field i.e.  $\| h \|^2 = -1$ on $\mathcal{Z}$.
In this way   $\pi': \mathcal{Z} \to  G_2^+(\R^4_1)$ may be viewed as  an analogous the usual Hopf
fibration  $\bb{S}^3 \to \bb{S}^2$. \\
 
 Next we define the  horizontal distribution $\goth{h} \subset T  \mathcal{Z}$ as the orthogonal complement with respect to the normal metric of the vertical distribution namely,
\begin{equation}\label{defHorizZ}
 \goth{h}_{(v,w)}= h_{(u,v)}^{\bot} \subset  T_{(v,w)} \mathcal{Z}, \,\, \forall (v,w) \in
 \mathcal{Z},
\end{equation}
thus
  $\goth{h} \subset T
\mathcal{Z}$ is  the complementary subbundle of the fibres of $\pi'$ and we decompose 
\begin{equation}\label{H}
 T_{(v,w)} \mathcal{Z} = \goth{h}_{(v,w)} \overset{\bot}{\oplus} \R (w,v),
 \end{equation}
 where $Ker \,\, d \pi'_{(v,w)} = \R h_{(v,w)} = \R (w,v)$.
    It is not difficult to verify that   $\goth{h}$
 defines a connection on 
   the principal bundle $ SO_o(1,1) \to \mathcal{Z} \to G_2^+
   (R^4_1)$. This is an example of sub-semi-riemannian geometry in which the metric on the distribution $\goth{h}$ is non-definite.

  \begin{lemma} \cite{huletttwistor}  For every  $(v,w) \in \mathcal{Z}$ the subspaces   $Ker\, d \pi_{(v,w)} $
 and $Ker \, d \pi'_{(v,w)} $ are orthogonal w.r.t. the normal metric  \eqref{inn.prod.so(3,1)}.
 \end{lemma}
  From the preceding Lemma and \eqref{H}  it  follows that $  Ker \, d \pi_{(v,w)} \subset
 \goth{h}_{(v,w)}$. Thus we have the   orthogonal decomposition
\begin{equation} \label{ortdecH}  \goth{h}_{(v,w)}=  Ker \, d \pi_{(v,w)} \overset{ \bot }{\oplus} L_{(v,w)}, 
\end{equation}
in which the subspace $L_{(v,w)}$ is the horizontal lift via $d \pi $ of the oriented spacelike $2$-plane $V = [v \wedge w]^{\bot} \subset T_v \bb{S}^3_1$.
 Note that  the
 normal metric~\eqref{inn.prod.so(3,1)} restricted to $\goth{h}$ has signature $(++)$ on $L$,
 and $(--)$ on $Ker \, d \pi$.  Hence the invariant  metric on $\goth{h}$ is neutral i.e. has
signature $(++--)$.\\
 
 The  geometry of the twistor space $\mathcal{Z} $ may be studied  with the aid of the so-called  Maurer-Cartan form $\beta$ of $\mathcal{Z}$  introduced by Burstall and Rawnsley in~\cite{burstall-rawnsley} of which we give a brief account. \\
Let $\goth{g} = \goth{so}(3,1)$ and recall   the
reductive decomposition $ \goth{g}= \goth{k} \oplus \goth{p}$. Consider the surjective application  $\xi_o : \goth{g}  \ni X \mapsto  \frac{d}{dt}|_{t=0} \exp(tX).o \in T_o \mathcal{Z}$. It follows that $\xi_o$ has  kernel  $\goth{k}$ and  restricts to an isomorphism $\goth{p} \to  T_o \mathcal{Z}$. 
Now form the associated vector bundle $ [\goth{p}] := SO_o(3,1) \times_K \goth{p}$. Then the map 
$$    [(g, X)] \mapsto \frac{d}{dt}|_{t=0} \exp(tAd(g)X)x =  d \tau_g (\frac{d}{dt}|_{t=0} \exp(tX).o ), \quad x= g.o,$$
 establishes an isomorphism of the associated bundle $[\goth{p}]$ and the tangent bundle $T \mathcal{Z}$, where $\tau_g$ is the isometry of $\mathcal{Z}$ sending $g'.o$ to $gg'.o$. 
 
Since $\goth{p}$ is an $Ad(K)$-invariant subspace of  $\goth{g}$,  one has the inclusion $[\goth{p}] \subset [\goth{g}]:= \mathcal{Z} \times \goth{g}$,  given by $ [\goth{p}] \ni [(g,X)] \mapsto (g.o, Ad(g)X) \in [\goth{g} ]$. 
Note that the fiber of $[\goth{p}] \to \mathcal{Z}$ over the point $g.o$ identifies with 
$ \{ g.o\} \times Ad(g)\goth{p} \subset [\goth{g}]$. This shows that  
 there exists an identification of  $T \mathcal{Z}$  with a subbundle of the trivial bundle $[\goth{g}]$. This inclusion   may be viewed as an $\goth{g} $-valued one-form  on $\mathcal{Z}$ which will be denoted by $\beta$.   
Note that every   $X \in \goth{g} $ determines a flow on $\mathcal{Z}$ defined by $\varphi_t (x) = \exp(tX).x$, which  is an isometry of $\mathcal{Z}$ for any $t \in \R$. The vector field of the flow is then a Killing field denoted by $X^*$ which is  given by 
$$ X^*_x =   \frac{d}{dt}|_{t=0}
\exp(tX) x, \quad \forall x \in \mathcal{Z}.    
$$   
It is not difficult to show that 
$$    \beta_x(X^*_x) = Ad(g)[Ad(g^{-1})(X)]_{\goth{p}}, \quad \forall X \in \goth{p},$$
 where $x=g.o \in \mathcal{Z}$. In particular at $o \in \mathcal{Z}$ we have $ \beta_o(X^*_o)= X$ for any $X \in \goth{p}$.  From this formula it follows the equivariance of $\beta$ which is expressed by 
 \begin{equation}\label{equivariantbeta}
    \beta \circ d \tau_g = Ad(g) \beta, \quad  \forall g \in SO_o(3,1).
    \end{equation}
For   $x = g.o \in \mathcal{Z}$  the application  
  $\xi_x :
 \goth{g}  \to T_{x} \mathcal{Z}$ such that 
$
        X \overset{\xi_x}{\longmapsto} X^*_x$, 
 maps $\goth{so}(3,1)$ onto $ T_{x} \mathcal{Z}$,    and  
  restricts to an isomorphism $ Ad(g)(\goth{p})   \to  T_x \mathcal{Z}$ whose  inverse  conicides with 
  $\beta_{x}$. Note that $\xi$ satisfies $d \tau_g \circ \xi_o (X) = \xi_{g.o} (Ad(g) X)$, which is equivalent to \eqref{equivariantbeta}. 
   More details and properties of the one form $\beta$ and  proofs can be found in~\cite{burstall-rawnsley}.\\

Now recall the definition of the horizontal distribution $\goth{h} \subset T \mathcal{Z}$ given in \eqref{defHorizZ}. At  the basepoint $o=(e_1, e_4) \in \mathcal{Z}$ the subspace  determined by the horizontal distribution identifies with      
\begin{equation} \label{Ho}
\mathcal{H} = \{     \begin{pmatrix}
  0 & x & y & 0 \\
  -x& 0 & 0 & z\\
  -y& 0 & 0 & w\\
  0 & z & w & 0\\
  \end{pmatrix}, \quad x,y,z,w \in \R \} \subset \goth{p},
  \end{equation}
which is an 
   $Ad(K)$-invariant subspace of $\goth{p}$.\\ 
 
  The one-form  $\beta$ transfers the metric on the
fibers of $T \mathcal{Z}$, the horizontal distribution $\goth{h}$  to the fibers of $[\goth{p}]$. 
From the definition of $\beta$ above we get   
\begin{equation} \label{betaH}
        \beta_{g.o}(\goth{h}_{g.o}) =   Ad(g) \mathcal{H} = [\mathcal{H}]_{g.o}, \quad \forall g \in SO_o(3,1).
        \end{equation}
        Thus $\beta$ identifies  the horizontal distribution $\goth{h} \subset T \mathcal{Z}$ with the sub bundle  $[\mathcal{H}] \subset [\goth{p}]$. \\
 
\subsection{ Invariant almost complex structures on $\goth{h}$.} 

  Let   $L_{(v,w)} \subset \goth{h}_{(v,w)}$ be the  $d \pi$-horizontal lift  of the oriented spacelike $2$-plane $V=[v \wedge w]^{\bot} \subset T_v \bb{S}^3_1$. 
 Denote by  $J^L_{(v,w)} :   L_{(v,w)} \to L_{(v,w)} $  the $ d \pi$-horizontal lift   of the positively oriented $\frac{\pi }{2}$-rotation 
 on the spacelike $2$-plane $V=[v \wedge w]^{\bot} $.\\ 
On the other hand let $J^{V}_{(v,w)}$ be the complex structure on the tangent space $Ker \, d \pi_{(v,w)}$ of the fibre of $v$ (recall that the fibers of $\pi : \mathcal{Z} \to \bb{S}^3_1$ are hyperbolic $2$-spaces, hence complex manifolds).  
Both structures together yield an  almost
complex structure $J'$  on the distribution $ \goth{h} = L \oplus Ker \, d \pi$ which is defined by
 \begin{equation}\label{thestructureJ}   J' = \left \{
       \begin{array}{ll}
        J^L ,& \text{on} \, L\\
        J^{\mathcal{V}},& \text{on}\, Ker \, d \pi \\
        \end{array}
        \right. 
        \end{equation}
  At  the base point $o= (e_1, e_4)$ it is possible to describe explicitly  the action of $J$. In fact  
$$
L_o = \{ \begin{pmatrix}  0 & x & y & 0\\
                          -x& 0 & 0 & 0\\
                          -y& 0 & 0 & 0\\
                          0 & 0 & 0 & 0 \\
                          \end{pmatrix} \}, \quad
                          Ker \, d\pi_o = \{\begin{pmatrix}
                           0 & 0 & 0 & 0\\
                          0& 0 & 0 & z\\
                          0& 0 & 0 & w\\
                          0 & z & w & 0 \\
                          \end{pmatrix} \}.
                          $$
On $L_{o}$ the complex structure is obtained by lifting via $d
\pi|_o$ the oriented rotation  on $[e_2 \wedge e_3]=[e_1 \wedge e_4]^{\bot} \subset
T_{e_1} \bb{S}^3_1$ given by  
$$ e_2 \mapsto e_3, \, e_3 \mapsto -e_2.$$

The complex structure $J^{\mathcal{V}}$ on the fibre  
$$\mathcal{Z}_{e_1} = \pi^{-1}(e_1) \equiv \bb{H}^2 \equiv SO(2,1)/SO(2) $$  is given at $o \in \mathcal{Z}_{e_1}$ by
$$
J^{\mathcal{V}}: \begin{pmatrix}           0 & 0 & 0 & 0\\   
                          0 & 0 & 0 & z\\
                          0&  0 & 0 & w\\
                          0& z & w & 0 \\
                          \end{pmatrix} \mapsto
                          \begin{pmatrix}
                          0 & 0 & 0 & 0\\
                          0& 0 & 0 & -w \\
                          0& 0 & 0 & z \\
                          0& -w & z & 0 \\
                          \end{pmatrix}.    $$
                       
                   According to the definition given before,  the almost complex structure $J' : \mathcal{H} \to  \mathcal{H}$ is given  by  
\begin{equation}\label{J1}
 J' \begin{pmatrix}
                           0 & x & y & 0\\
                          -x& 0 & 0 & z\\
                          -y& 0 & 0 & w\\
                          0 & z & w & 0 \\
                          \end{pmatrix} =
                          \begin{pmatrix}
                          0 & -y & x & 0\\
                          y& 0 & 0 & -w \\
                          -x& 0 & 0 & z \\
                          0 & -w & z & 0 \\
                          \end{pmatrix}.
                          \end{equation}
 
 The  second  almost complex structure $J''$ on $\goth{h}$ is obtained by reversing the complex structure on the fibers, namely
 \begin{equation}\label{thestructureJ'}   J'' = \left \{
       \begin{array}{ll}
        J^L ,& \text{on} \, L\\
        -J^{\mathcal{V}},& \text{on}\, Ker \, d \pi \\
        \end{array}
        \right. 
        \end{equation}
 The action of $J''$ on the subspace $\mathcal{H} \subset \goth{p}$ is given by 
\begin{equation}\label{J''}
 J'' \begin{pmatrix}
                           0 & x & y & 0\\
                          -x& 0 & 0 & z\\
                          -y& 0 & 0 & w\\
                          0 & z & w & 0 \\
                          \end{pmatrix} =
                          \begin{pmatrix}
                          0 & -y & x & 0\\
                          y& 0 & 0 & w \\
                          -x& 0 & 0 & -z \\
                          0 & w & -z & 0 \\
                          \end{pmatrix}.
                          \end{equation} 
 
We summarize our discussion above in the following   
 \begin{lemma}   Let $J$ be either $J'$ or $ J''$. Then $J$  commute with 
  $ \{ Ad(x)|_{\mathcal{H}} : x \in K \}$  and  is  orthogonal i.e.  
$$
      \la JX, JY \ra = \la X,  Y \ra, \quad \forall X,Y \in \mathcal{H},  $$
      where  $J \in \{ J', J'' \}$. 
      Thus $J'$ and $J''$ are $SO_o(3,1)$-invariant almost complex structures on  $ \goth{h}= [ \mathcal{H} ] \subset T \mathcal{Z}$.
      \end{lemma} 
As a consequence of a theorem by LeBrun  it can be shown that 
the almost CR structure $(  \goth{h}= [ \mathcal{H} ], J')$ on $\mathcal{Z}$ is integrable (see~\cite{kobak}).

\section{ Harmonic maps into $\mathcal{Z}$ }  

Here we study the harmonic map equation for smooth maps
 $ \phi :M \to \mathcal{Z}$  from a Riemann surface into the twistor bundle, on which  we have fixed  the normal metric $\la .,. \ra$.  
 Let  $\Omega \subset M$ be a relatively compact domain of $M$ and define the twistor energy of $\phi$ over $\Omega$ by
\begin{equation}\label{energyonomega}
E_{\Omega} (\phi) = \frac{1}{2}\int_{\Omega} \| d \phi \|^2 dA,
\end{equation}
   where $dA$ is the area form  on $M$ determined by a  conformal metric $g$, and  $\| d \phi \|^2$ is the Hilbert-Schmidt norm of $d \phi$ defined by  $\| d \phi \|^2 = \sum_i \la d \phi (e_i), d \phi (e_i) \ra$ for any orthonormal frame $\{ e_i \}$ on $M$.
 By definition $\phi$ is harmonic if it is an extreme of the energy functional $ \phi \mapsto     E_{\Omega} (\phi)$  over  all relatively compact subdomains $\Omega \subset M$.    
 Note that the energy may be negative since the metric on $\mathcal{Z}$ is indefinite.\\
 
Let us denote by  $\nabla$ the levi-Civita connection on $  \mathcal{Z}$ determined by the normal metric, and  by $\nabla^{\phi}$ the induced connection on the pull-back bundle $\phi^{-1} T \mathcal{Z}$.  By the formula of the first variation of the energy ( see \cite{eells-lemaire}) it follows that $\phi$ is harmonic if and only if  its tension vanishes: $tr(\nabla d\phi) =0$. 
If $M$ is  Riemann surface then     $\phi$ is harmonic  if and only if  on every  local complex coordinate $z$ on $M$ the following equation holds
\begin{equation} \label{generalharmoniceq}
        \nabla^{\phi}_{\Bz} d \phi(\Pz) = 0,        
\end{equation}
where the left hand of this equation is just a non-zero multiple of the  tension field of  $\phi$.   \\
 
  For our purposes we need a reformulation  of equation~\eqref{generalharmoniceq}     in terms of  the Maurer-Cartan form or Moment map $\beta$ of $\mathcal{Z}$, see~\cite{black},~\cite{burstall-rawnsley}. \\   
   Let  $D$ be the {\it canonical connection of second kind} i.e. the affine connection  on $\mathcal{Z}$ for which the parallel transport along the curve  $t \to \exp(tX).x$ is realized  by $ d \exp(tX) $. Hence at the basepoint $o \in \mathcal{Z}$ we have 
   \begin{equation} \label{canonicconnectionD}     D_{X^*} Y^*(o) = \frac{d}{dt}|_{t=0}  d \exp(-t X) Y^* = [X^*, Y^*](o) = -[X,Y]_{\goth{p}}. 
   \end{equation}
   
   Note also that $D$ is determined by  the condition  $(D_{X^*}  X^*)_o =0, \forall X \in \goth{p}$.  Since $\nabla \neq D$ and  $D \la \,, \, \ra =0$, $D$  has non-vanishing torsion. From\eqref{canonicconnectionD} we obtain   
  \begin{equation}\label{torsionD}
  T_o^D(X^*,Y^*) = - [X,Y]_{\goth{p}}, \quad X,Y \in \goth{p}.
  \end{equation}
 
  The following formula allows to compute  $D$ in terms of $\beta $ and the Lie algebra structure of $\goth{so}(3,1)$, 
  
  \begin{lemma}~\cite{burstall-rawnsley}   
  \begin{equation} \label{betaformula}     \beta(D_X Y) = X \beta(Y) - [ \beta(X), \beta(Y)], \quad X,Y \in \mathcal{X} (\mathcal{Z}). 
  \end{equation}
  \end{lemma}
\bigskip

Let  us now compute the Levi-Civita connection $\nabla$ of  the normal metric on $\mathcal{Z}$. Since $P: SO_o(3,1) \to \mathcal{Z}$ is a (pseudo) riemannian submersion, $P$ sends $\goth{p}$-horizontal geodesics in $SO_o(3,1)$ onto $\nabla$-geodesics in $\mathcal{Z}$. Hence $\nabla_{X^*}X^* =0$ for every $X \in \goth{p}$ which implies  $\nabla_{X^*}Y^* + \nabla_{Y^*}X^* =0$, for any $X,Y \in \goth{p}$. Now since $\nabla$  is 
torsionless we get  
\begin{equation}\label{nablanormal}
 \nabla_{X^*}Y^*= \frac{1}{2} [X^*, Y^*], \forall X,Y \in \goth{p}.
 \end{equation}
 Hence at $o \in \mathcal{Z}$ we get
 \begin{equation}\label{nablanormalato}
 (\nabla_{X^*}Y^*)(o)= -\frac{1}{2} [X,Y]_{\goth{p}}, \,\, \forall X,Y \in \goth{p}.
 \end{equation}

   We are now ready to obtain a  formula for the Levi-Civita connection $\nabla$ on $\mathcal{Z}$ in terms of $\beta$ namely,
\begin{lemma}
   \begin{equation}\label{formulaLeviCivita}
   \beta (\nabla_X Y) =  X \beta(Y) - [ \beta(X), \beta(Y)] + \frac{1}{2} \pi_{\goth{p}} ([ \beta(X), \beta(Y)]), \quad X,Y \in \mathcal{X}(\mathcal{Z}),
   \end{equation}
   where $\pi_{\goth{p}} : \mathcal{Z} \times \goth{so}(3,1) \to SO_o (3,1) \times_K \goth{p}\equiv T\mathcal{Z}$ is the projection onto the  the tangent bundle of $\mathcal{Z}$.
   \end{lemma}
    \noindent {\bf Proof.} Let $X^*, Y^* $ be  Killing vector fields on $\mathcal{Z}$ determined by $X,Y \in \goth{p}$. \\
From de definition of $\beta$,  \eqref{canonicconnectionD} and \eqref{nablanormalato} we have 
    $$    
   \beta ((\nabla_{X^*} Y^*)(o)) - \beta ((D_{X^*} Y^*)(o)) = -\frac{1}{2} [X,Y]_{\goth{p}}+ [X,Y]_{\goth{p}}=
   \frac{1}{2} [X,Y]_{\goth{p}}, \,\,  \forall X,Y \in \goth{p}. $$ 
On the other hand the difference tensor $\nabla -D $ is  $SO_o(3,1)$-invariant, and so is  $\beta(\nabla-D) = \beta(\nabla) - \beta(D)$ by formula \eqref{equivariantbeta}. Hence it is determined by its value at the point $o \in \mathcal{Z}$. Thus  formula  \eqref{formulaLeviCivita} follows. $\square$ \\ 
  
Define the $D$-fundamental form of $\phi :M \to \mathcal{Z}$ by 
    \begin{equation} \label{D2ndfundform}
               D d\phi (U,V) =  D^{\phi}_{U} d \phi(V) - d \phi (\nabla^M_U V), \quad  U,V \in \mathcal{X}(M),
    \end{equation}
   in which $D^{\phi}$ is the connection on $\phi^{-1} T \mathcal{Z}$ determined by $D$  and $\nabla^M$ is the Levi-Civita connection on $M$ determined by a conformal metric.  The map $\phi :M \to \mathcal{Z}$ is called $D$-harmonic if and only if   $tr (D d\phi) =0$ or equivalently 
   $$D^{\phi}_{\Bz} d \phi(\Pz)=0.$$ 
   
 \begin{lemma}\label{harmonicmaplemma}     $\phi :M \to \mathcal{Z}$ is harmonic if and only if   
 \begin{equation} \label{Zharmonicmapequation}
 \bz (\phi^* \beta)' -
 [(\phi^* \beta)''  \wedge (\phi^* \beta)' ] =0,
 \end{equation}
 where $\phi^* \beta = (\phi^* \beta)' + (\phi^* \beta)''$ is the decomposition of the (complex) one form $\phi^* \beta$ into one forms of type $(1,0)$ and $(0,1)$.   
 \end{lemma}
 \noindent{\bf Proof.} It is consequence of the following  formula for the tension of $\phi$ 
  \begin{equation} \label{betatorsionphi}
  \begin{array}{cc}
 -\beta ( \nabla_{\Bz} d \phi (\Pz) ) dz \wedge d \bar{z} = \bz (\phi^* \beta)' -
 [(\phi^* \beta)''  \wedge (\phi^* \beta)' ].\\
 \end{array}
 \end{equation}
 To obtain formula \eqref{betatorsionphi} we note first  that from  \eqref{formulaLeviCivita} we have $\beta (tr \nabla d \phi) =  \beta (tr D d \phi)$ which since $M$ is a Riemann surface is equivalent to 
 $$      \beta(\nabla^{\phi}_{\Bz} d \phi (\Pz) ) = \beta(D^{\phi}_{\Bz} d \phi (\Pz)).  $$     
 On the other hand  from formula \eqref{betaformula} we obtain   
 \begin{equation} \label{betatorsionphi2}
 \beta(D^{\phi}_{\Bz} d \phi (\Pz)) = \Bz \beta d \phi(\Pz) - [ \beta d \phi(\Bz), \beta d \phi(\Pz) ],      
 \end{equation}   
 from which \eqref{betatorsionphi} follows. $\square$ 

\section{Harmonic twistor lifts  }
 A natural question  is to characterize those conformaly immersed surfaces  $f: M  \to \bb{S}^3_1$   whose  twistor lift $\widehat{f} :M \to \mathcal{Z}$ is a harmonic map.   \\
Let $f: M \to \bb{S}^3_1$ be a conformal immersed surface and 
  $F \in  SO_o (3,1)$ be a frame of $f$ which is adapted to the surface $f$ i.e. 
 $F$ is  a locally defined map  on an open subset  $U \subset M$ satisfying 
      \begin{equation*}
    \begin{array}{ll}
    f(x)= F_1(x),\,\, F_4(x)=  \n (x), \\ 
     \text{span} \{ F_2(x), F_3(x) \} = df_x (T_xM), \forall  x \in U,\\
     \end{array}
     \end{equation*} 
     where  $F_i = F. e_i$ are the columns of the matrix $F$.
   From  the  structure equations \eqref{gauss-codazzi} of the immersed spacelike surface $f$  we obtain the following 	evolution equations of the frame $F$   
\begin{equation}\label{structure eqs}
\begin{array}{l}
f_z = \frac{e^u}{\sqrt{2}}\, (F_2-i F_3 ),\\
\\
  (F_2)_z=
-\frac{e^{u}}{\sqrt{2}}.f -i u_z F_3 + ( \frac{ e^{-u} \xi+
e^{u}H }{\sqrt{2}}).\n
\\
\\
 (F_3)_z= i \frac{e^{u}}{\sqrt{2}}.f + i  u_z F_2 + i  ( \frac{
e^{-u} \xi- e^{u}H }{\sqrt{2}}).\n \\
\\
\n_z = (\frac{e^{-u} \xi+ e^u H }{\sqrt{2}}). F_2 +i (\frac{e^{-u} \xi-e^u H }{\sqrt{2}}).F_3\\
\\
\end{array}
\end{equation}

Written in matrix form these equations take the form 
\begin{equation} \label{FAB}
   F_z = F A, \quad  F_{\bar{z}} = F B,
\end{equation} 
where the complex matrices $A,B= \overline{A}$ are given by     
\begin{equation}\label{matrA}
A =
\begin{pmatrix}
0   &   -\frac{e^u}{\sqrt{2}} & i\frac{e^u}{\sqrt{2}} & 0 \\
\frac{e^u}{\sqrt{2}}& 0 &  i u_z & \frac{ e^{-u} \xi + e^u H}{\sqrt{2}}\\
-i \frac{e^u}{\sqrt{2}} & -i u_z & 0 & i (  \frac{e^{-u} \xi-e^u H}{\sqrt{2}})\\
0 &   \frac{ e^{-u} \xi+ e^{u}H }{\sqrt{2}} &  i ( \frac{e^{-u}
\xi-e^u H}{\sqrt{2}})&
0 \\
\end{pmatrix}
\end{equation}

\begin{equation}\label{matrB}
B=
\begin{pmatrix}
0   &   -\frac{e^u}{\sqrt{2}} & -i\frac{e^u}{\sqrt{2}} & 0 \\
\frac{e^u}{\sqrt{2}}& 0 &  -i u_{\bar{z}} & \frac{ e^{-u} \overline{\xi} + e^u H}{\sqrt{2}}\\
i \frac{e^u}{\sqrt{2}} & i  u_{\bar{z}} & 0 & -i (  \frac{e^{-u} \overline{\xi}-e^u H}{\sqrt{2}})\\
0 &   \frac{ e^{-u} \overline{\xi}+ e^{u}H }{\sqrt{2}} &  -i (
\frac{e^{-u} \overline{ \xi}-e^u H}{\sqrt{2}})&
0 \\
\end{pmatrix}.
\end{equation}
The pull-back $\alpha = F^{-1} dF$ of the Maurer-Cartan form of $SO_o(3,1)$  is given in terms of $A,B$ by 
\begin{equation}\label{pullbackAB}
\alpha= A dz + B d\bar{z}.
\end{equation}Taking $\goth{p}$ and $\goth{k}$ projections  we obtain $\alpha= \alpha_{\goth{k}} + \alpha_{\goth{p}}$, with   
\begin{equation} \label{alphakp}
\alpha_{\goth{k}} = A_{\goth{k}} dz + B_{\goth{k}} d\bar{z}, \quad  
\alpha_{\goth{p}} = A_{\goth{p}} dz + B_{\goth{p}} d\bar{z},
\end{equation}
in which  
\begin{equation}\label{dzlift}
   A_{\goth{p}}  = 
   \begin{pmatrix}
0   &   -\frac{e^u}{\sqrt{2}} & i\frac{e^u}{\sqrt{2}} & 0 \\
\frac{e^u}{\sqrt{2}}& 0 &  0 & \frac{ e^{-u} \xi + e^u H}{\sqrt{2}}\\
-i \frac{e^u}{\sqrt{2}} & 0 & 0 & i (  \frac{e^{-u} \xi-e^u H}{\sqrt{2}})\\
0 &   \frac{ e^{-u} \xi+ e^{u}H }{\sqrt{2}} &  i ( \frac{e^{-u}
\xi-e^u H}{\sqrt{2}})&
0 \\
\end{pmatrix},
\end{equation}
\begin{equation}\label{dzbarlift}
   B_{\goth{p}} =
   \begin{pmatrix}
0   &   -\frac{e^u}{\sqrt{2}} & -i\frac{e^u}{\sqrt{2}} & 0 \\
\frac{e^u}{\sqrt{2}}& 0 &  0 & \frac{ e^{-u} \overline{\xi} + e^u H}{\sqrt{2}}\\
i \frac{e^u}{\sqrt{2}} & 0 & 0 & -i (  \frac{e^{-u} \overline{\xi}-e^u H}{\sqrt{2}})\\
0 &   \frac{ e^{-u} \overline{\xi}+ e^{u}H }{\sqrt{2}} &  -i (
\frac{e^{-u} \overline{ \xi}-e^u H}{\sqrt{2}})&
0 \\
\end{pmatrix},
\end{equation}
\begin{equation} \label{AkBk}
  A_{\goth{k}} = 
  \begin{pmatrix}
0   &   0 & 0 & 0 \\
0& 0 &  i u_z & 0\\
0 & -i u_z & 0 & 0\\
0 &   0 & 0 & 0 \\
\end{pmatrix}, \quad 
B_{\goth{k}}=
\begin{pmatrix}
0   &   0 & 0 & 0 \\
0& 0 &  -i u_{\bar{z}} & 0\\
0 & i  u_{\bar{z}} & 0 & 0\\
0 &   0 & 0 & 0 \\
\end{pmatrix}.
\end{equation} 

Also since $M$ is a Riemann surface we decompose $\alpha_{\goth{p}}$ and  $\alpha_{\goth{k}}$ into its $(1,0)$ and $(0,1)$ parts 
 $$ \alpha_{\goth{p}} = \alpha'_{\goth{p}} + \alpha''_{\goth{p}}, \quad 
 \alpha_{\goth{k}} = \alpha'_{\goth{k}} + \alpha''_{\goth{k}}.$$
 
 In terms of the above matrices we have
\begin{equation}
\alpha'_{\goth{k}} = A_{\goth{k}} dz, \,  \alpha''_{\goth{k}}= B_{\goth{k}} d\bar{z}, \quad 
\alpha'_{\goth{p}}= A_{\goth{p}} dz,  \,   \alpha''_{\goth{p}} = B_{\goth{p}} d\bar{z}.
\end{equation}
Note that  the adapted frame $F$ of $f$ is also a frame of the twistor lift $\widehat{f}$ since  $F.o = (F_1, F_4) = (f, \n) = \widehat{f}$, where $o=(e_1, e_4) \in \mathcal{Z}$ is the fixed basepoint.   
 From  formula $\widehat{f}^* \beta = Ad(F) \alpha_{\goth{p}}$  we obtain $$(\widehat{f}^* \beta)' = Ad(F) \alpha'_{\goth{p}}, \quad (\widehat{f}^* \beta)'' = Ad(F) \alpha''_{\goth{p}}.$$
Now using the identity (see \cite{burstall-pedit} pag. 241)
$$  d Ad(F) = Ad(F) \circ ad \, \alpha , $$
  we compute 
 \begin{equation}
 \begin{array}{cc}
      \bz (\widehat{f}^* \beta)' = \bz \{  Ad(F) \alpha'_{\goth{p}}    \} =  Ad(F) \{   \bz \alpha'_{\goth{p}} + [\alpha_{\goth{k}} \wedge \alpha'_{\goth{p}} ] + [\alpha''_{\goth{p}} \wedge \alpha'_{\goth{p}} ]   \}.\\
      \end{array}
      \end{equation}

On the other hand 
 \begin{equation}
 \begin{array}{cc} 
  [ (\widehat{f}^* \beta)'' \wedge (\widehat{f}^* \beta)'] = Ad(F)   [\alpha''_{\goth{p}} \wedge \alpha'_{\goth{p}} ] \\
 \end{array}
 \end{equation} 
 
 
Therefore  
$$  \bz (\widehat{f}^* \beta)' - [ (\widehat{f}^* \beta)'' \wedge (\widehat{f}^* \beta)'] = Ad(F) \{\bz \alpha'_{\goth{p}} + [\alpha_{\goth{k}} \wedge \alpha'_{\goth{p}} ]     \}
$$    
Thus a direct consequence of Lemma 5  is 
\begin{proposition}
Let $f:M \to \bb{S}^3_1$ be a conformal immersed surface and let $\widehat{f}:M \to \mathcal{Z}$ be its twistor lift. Then $\widehat{f}$ is a harmonic map if and only if 
for every adapted frame $F$ of $f$ the $\goth{so}(3,1)$-valued one form $\alpha = F^{-1}dF$ satisfies 
\begin{equation}\label{harmeqalpha}
        \bz \alpha'_{\goth{p}} + [\alpha_{\goth{k}} \wedge \alpha'_{\goth{p}} ] =0.
\end{equation}
\end{proposition}

 \medskip
 
A  reformulation of \eqref{harmeqalpha}  in terms of the evolution matrices $A,B$ is easily obtained. Recall that  $ \alpha'_{\goth{p}} = A_{\goth{p}} dz$ and  $\alpha''_{\goth{p}} = B_{\goth{p}} d \bar{z}$, hence 
 
 \begin{equation} \label{auxbeta}
 \bz \alpha'_{\goth{p}} + [\alpha_{\goth{k}} \wedge \alpha'_{\goth{p}} ] = 
 - ( \Bz A_{\goth{p}} + [ B_{\goth{k}} ,  A_{\goth{p}}] \,) \,   dz  \wedge d \bar{z} 
 \end{equation}

\begin{corollary}
Let $f: M \to \bb{S}^3_1$ be a conformal immersion and let $\widehat{f}:M \to \mathcal{Z}$ be the twistor lift of $f$. Then  $\widehat{f}$ is harmonic if and only if for every adapted frame $F$ of $f$ the complex matrices $A,B$ defined by $\alpha= F^{-1} dF = A dz + B d \bar{z}$ satisfy 
\begin{equation}\label{harm.map.eq.fhat}
    \Bz  A_{\goth{p}} + [ B_{\goth{k}}, A_{\goth{p}}] =0,
\end{equation}
\end{corollary}

From the  explicit form of the matrices $A_{\goth{p}}$ and $B_{\goth{k}}$ in \eqref{dzlift} and \eqref{AkBk} we conclude from equation \eqref{harm.map.eq.fhat} that  $\widehat{f}$ is harmonic if and only if  
\begin{equation} 
\begin{array}{cc}
\Bz ( e^{-u} \xi + e^{u} H ) +  u_{\bar{z}} (e^{-u} \xi - e^{u} H)=0,\\
\Bz (e^{-u} \xi - e^{u} H) + u_{\bar{z}}  ( e^{-u} \xi + e^{u} H )=0.
\end{array}
\end{equation} 
Cancelling terms we are left with 
\begin{equation}\label{systemHE}
\begin{array} {ll}
 e^{-u}  \xi_{\bar{z}} + e^{u}  H_{\bar{z}} =0, \\
 e^{-u}  \xi_{\bar{z}} -  e^{u}  H_{\bar{z}} =0. \\
 \end{array}
 \end{equation}
   Combining with  Codazzi's equation $\xi_{\bar{z}}= e^{2u} H_z$,   system \eqref{systemHE} is equivalent to    
 \begin{equation}
 \begin{array}{ll}
  e^{u} (H_z + H_{\bar{z}}) =0, \\
  e^{u} (H_z - H_{\bar{z}}) =0,\\
  \end{array}
 \end{equation}
hence  $H_x = H_y =0$ and conversely.  We have thus proved the following 
   \begin{theorem}\label{harmonictwsitorlift}
  Let $f: M \to \bb{S}^3_1$ be a conformal (hence spacelike) immersion and let $\widehat{f}:M \to \mathcal{Z}$ be its twistor lift. 
  Let $\la . , . \ra$ be the normal metric  on the twistor space $\mathcal{Z}$.  
  Then $\widehat{f}: M \to (\mathcal{Z} ,\la . , . \ra )$ is a harmonic map if and only if the immersed  surface $f$ has constant mean curvature.
  \end{theorem}
   By Codazzi's equation $\widehat{f}$ is harmonic if and only if  the Hopf complex differential   $q= \xi dz^2= - \la
 f_{z z}, \n \ra^{c} dz \otimes dz$ is holomorphic.

\section{One parameter deformations}

Here we show that harmonic twistor lifts exist within a family parameterized by the complex numbers of unit modulus. \\
 
Let $f: M \to \bb{S}^3_1$ a conformal immersion and $F$  a (local) adapted frame of $f$, hence $F$ is also a frame of $\widehat{f}$. Let  $\alpha = F^{-1}dF$ be the pullback of the Maurer-Cartan form by $F$.  According to the reductive decomposition  $\goth{so}(3,1) = \goth{k} \oplus \goth{p}$ we 
 decompose as before $\alpha =  \alpha_{\goth{k}} + \alpha_{\goth{p}}$ where in terms of   matrices $A_{\goth{p}}, B_{\goth{p}}$, $A_{\goth{k}}, B_{\goth{k}}$~\eqref{dzlift},~ \eqref{dzbarlift}  and~\eqref{AkBk} these forms are expressed by 
 \begin{equation*} 
 \begin{array}{ll}
  \alpha_{\goth{k}}= A_{\goth{k}} dz + B_{\goth{k}} d\bar{z}, & 
 \alpha_{\goth{p}}= A_{\goth{p}} dz + B_{\goth{p}} d\bar{z},\\
 \alpha'_{\goth{p}}=A_{\goth{p}} dz, & \alpha''_{\goth{p}}=B_{\goth{p}} d\bar{z}.
 \end{array}
 \end{equation*}   
On the other hand $\alpha$  satisfies the Maurer-Cartan equation  $  d \alpha + \frac{1}{2} [\alpha \wedge \alpha]=0$ which  splits up into
\begin{equation}\label{decompalpha}
\begin{array}{cc}     
    \bz \alpha'_{\goth{p}} + [ \alpha_{\goth{k}} \wedge \alpha'_{\goth{p}}] +
    \pz \alpha''_{\goth{p}} + [ \alpha_{\goth{k}} \wedge \alpha''_{\goth{p}}]+[ \alpha'_{\goth{p}} \wedge \alpha''_{\goth{p}}]_{\goth{p}}=0,\\
    d\alpha_{\goth{k}} + \frac{1}{2}[ \alpha_{\goth{k}} \wedge \alpha_{\goth{k}}]+
    [ \alpha'_{\goth{p}} \wedge \alpha''_{\goth{p}}]_{\goth{k}}=0.\\    
\end{array}
\end{equation}  
\medskip

We give below a further  property of the geometry  of our twistor lifts $\widehat{f}$ which is  consequence of  the form of the structure equations of the immersion $f$ which  is reflected in  the matrices $A_{\goth{p}}$~\eqref{dzlift} and   $\overline{ A_{\goth{p}}} = B_{\goth{p}}$~\eqref{dzbarlift}. 
 
\begin{lemma}\label{ontheconditionX}
Let $f: M \to \bb{S}^3_1$ be a conformal immersion,  $F$ an arbitrary  adapted frame of $f$, and  $z=x+iy$ a local complex coordinate on $M$. Set $F^{-1}F_z = A_{\goth{p}} + A_{\goth{k}}$ and $F^{-1}F_{\bar{z}} = B_{\goth{p}} + B_{\goth{k}}$, where the complex matrices $A_{\goth{p}}$, $B_{\goth{p}}$ given by~\eqref{dzlift},~ \eqref{dzbarlift}.   
 Then the one forms  $\alpha'_{\goth{p}} =A_{\goth{p}}dz$ and $ \alpha''_{\goth{p}}= B_{\goth{p}}d\bar{z}$   satisfy 
   \begin{equation} \label{conditionX}
    [\alpha'_{\goth{p}} \wedge \alpha''_{\goth{p}}]_{\goth{p}}=0.
    \end{equation}
\end{lemma} 
\noindent {\bf Proof.} 
 From the structure equations~\eqref{structure eqs} of the immersion $f$  we obtain $\la f_z, \n_{\bar{z}} \ra^c = \la f_{\bar{z}}, \n_{z} \ra^c$ and  $\la \Pz F_i - [\Pz F_i]^T, \n_{\bar{z}} \ra^c =0, \, i= 2,3$, 
 where $[\Pz F_i]^T$ denotes the projection onto the tangent bundle of the immersed surface. These equations clearly imply  $[A_{\goth{p}}, B_{\goth{p}}]_{\goth{p}}=0$. Hence 
$$   
[\alpha'_{\goth{p}} \wedge \alpha''_{\goth{p}}]_{\goth{p}}=[A_{\goth{p}}, B_{\goth{p}}]_{\goth{p}} \, dz \wedge d\bar{z}=0. \quad \square $$

 Assume now that  $f: M \to \bb{S}^3_1$ has constant mean curvature. Thus  $\widehat{f}$ is  a harmonic map by  Theorem~\ref{harmonictwsitorlift},  and so it satisfies the harmonic map equation  
$$ 0= \bz \alpha'_{\goth{p}} + [\alpha_{\goth{k}} \wedge \alpha'_{\goth{p}} ] = 
 - ( \Bz A_{\goth{p}} + [ B_{\goth{k}} ,  A_{\goth{p}}] \,) \,   dz  \wedge d \bar{z}.  $$        
Taking into account condition~\eqref{conditionX}  the first equation in~\eqref{decompalpha}  reduces to   
$$\pz \alpha''_{\goth{p}} + [ \alpha_{\goth{k}} \wedge \alpha''_{\goth{p}}]=0.$$  
Hence the pair of equations~\eqref{decompalpha} become 
\begin{equation*}
\begin{array}{ll}
(\textbf{a}) &\pz \alpha''_{\goth{p}} + [ \alpha_{\goth{k}} \wedge \alpha''_{\goth{p}}]=0,\\
(\textbf{b}) & d\alpha_{\goth{k}} + \frac{1}{2}[ \alpha_{\goth{k}} \wedge \alpha_{\goth{k}}]+
    [ \alpha'_{\goth{p}} \wedge \alpha''_{\goth{p}}]=0.\\
\end{array}
\end{equation*}

For  $\lambda \in \bb{C}$ with $|\lambda| =1$   set 
\begin{equation}\label{alphalambdaoneform} 
\lambda. \alpha = \alpha_{\lambda} = \lambda^{-1} \alpha'_{\goth{p}} + \alpha_{\goth{k}} + \lambda \alpha''_{\goth{p}}.
\end{equation}
Due to $\overline{A_{\goth{p}}}=B_{\goth{p}}$ and  $\overline{A_{\goth{k}}}=B_{\goth{k}}$, $\alpha_{\lambda}$ is $\goth{so}(3,1)$-valued for every $\lambda \in \bb{S}^1$.  
Moreover  $\lambda . \alpha = \alpha_{\lambda}$  defines an action of $\bb{S}^1$ on $\goth{so}(3,1)$-valued $1$-forms which leaves invariant the solution set of equations (\textbf{a}) and (\textbf{b}) above. 
Comparing coefficients of $\lambda $  it follows that  equations $(\bf a)$ and $(\bf b)$ above hold for $\alpha$ if and only if $\alpha_{\lambda}$ satisfies       
   $$  
  d \alpha_{\lambda} + \frac{1}{2} [\alpha_{\lambda} \wedge \alpha_{\lambda}]=0, \forall \lambda \in \bb{S}^1.  $$ 
  This is the so-called {\it zero curvature condition} (ZCC)~\cite{burstall-pedit}. 
In this way the harmonic map equation for twistor lifts to $\mathcal{Z}$  is encoded in a loop of "zero curvature" equations. \\

Now  let us assume  that the Riemann surface $M$ is simply connected (otherwise we  pass to  its  universal covering space $\tilde{M}$), and fix  a base point $m_o \in M$. Then for each $\lambda \in \bb{S}^1$ we can integrate and solve  
\begin{equation}\label{extendedframeeqs}
dF_{\lambda} = F_{\lambda} \alpha_{\lambda}, \quad F_{\lambda} (m_o)= Id.
\end{equation}
 
The solution map    
$F_{\lambda}=(f_{\lambda}, (F_{\lambda})_2, (F_{\lambda})_3, \n_{\lambda}) : M \to SO_o(3,1)$  is  called an {\it extended frame} and satisfies   
\begin{equation}\label{lambdamat}
F^{-1}_{\lambda} (F_{\lambda})_z = \lambda^{-1} A_{\goth{p}}+ A_{\goth{k}} , \quad F^{-1}_{\lambda} (F_{\lambda})_{\bar{z}} = \lambda B_{\goth{p}} + B_{\goth{k}}, \quad \forall \lambda \in \bb{S}^1. 
\end{equation}
Moreover since   $ (\alpha_{\lambda})_{\goth{p}} = \lambda^{-1} \alpha'_{\goth{p}} + \lambda \alpha''_{\goth{p}} =(\alpha_{\lambda})'_{\goth{p}} +(\alpha_{\lambda})''_{\goth{p}}$, and $(\alpha_{\lambda})_{\goth{k}}= \alpha_{\goth{k}}$,   
  then    the one form $\alpha_{\lambda}$ satisfies equations $(\textbf{a})$ and $(\textbf{b})$ for every $\lambda \in \bb{S}^1$.  Thus if     $P : SO_o(3,1) \to \mathcal{Z}$ denotes the projection map  $P (g) = g.o$, then    $\phi_{\lambda} =  P \circ F_{\lambda} :M \to \mathcal{Z}$ is harmonic  $\forall \lambda \in \bb{S}^1$.\\  
  The family  $\{ \phi_{\lambda}, \lambda \in \bb{S}^1 \}$ is called {\it the associated family} of the harmonic map $\widehat{f}$~\cite{burstall-rawnsley}.
  Note that   $\phi_{ \{ \lambda =1 \}} = \widehat{f} $, hence  each $ \phi_{\lambda} $ is a  deformation of $\widehat{f}$.  \\


Let  $f_{\lambda} = \pi \circ \phi_{\lambda} : M \to \bb{S}^3_1$, hence      from~\eqref{lambdamat}  we extract 
\begin{equation} \label{lambdafn}
\begin{array}{cc}
 (f_{\lambda})_z = \lambda^{-1} \frac{e^{u}}{\sqrt{2}} \left[ (F_{\lambda})_2-i(F_{\lambda})_3   \right],\\
 \\
 (\n_{\lambda})_z = \lambda^{-1} \left[ (\frac{e^{-u}\xi + e^{u}H}{\sqrt{2}})  (F_{\lambda})_2+i (\frac{e^{-u}\xi - e^{u}H}{\sqrt{2}}) (F_{\lambda})_3 \right].\\
 \end{array} 
\end{equation}

 From the first equation above we get 
 $$\la (f_{\lambda})_z ,(f_{\lambda})_z \ra =\la (f_{\lambda})_z ,(f_{\lambda})_{\bar{z}} \ra^c  = e^{2u}, $$
 
thus   $\{ f_{\lambda}, \lambda \in \bb{S}^1 \}$  is a family of conformal immersions into $\bb{S}^3_1$, with a common conformal factor $u$, hence  all $f_{\lambda}$ induce the same metric for every $\lambda \in \bb{S}^{1}$.  
Let  $H_{\lambda}$ be the mean curvature of $f_{\lambda}$. Since $u$ is the conformal parameter of $f_{\lambda}$, we get from~\eqref{lambdamat}, 
$$H_{\lambda} = -e^{-2u} \la (f_{\lambda})_{\bar{z}z}, \n_{\lambda} \ra =e^{-2u} 
 \la (f_{\lambda})_{z}, (\n_{\lambda})_z \ra = H, \,\,\,\,  \forall \lambda \in \bb{S}^1.
$$
Thus $H_{\lambda}$ does not depend on $\lambda$.
 On the other hand let $q_{\lambda} = \xi_{\lambda} dz \otimes dz$ be the Hopf complex differential of $f_{\lambda}$. Then from~\eqref{lambdamat} and~\eqref{lambdafn} we get 
\begin{equation}\label{hopflambda}
      \xi_{\lambda} =-\la (f_{\lambda})_{z z}, \n_{\lambda} \ra^c=  \la (f_{\lambda})_{z}, (\n_{\lambda})_z \ra^c = \lambda^{-2} \xi. 
\end{equation}
 Gauss equation for $f_{\lambda}$ reads 
\begin{equation}\label{gausscurvaturelambda}
K_{\lambda}= 1-H^2 + | \xi_{\lambda} |^2 e^{-4u}=  1- H^2 + |\xi |^2 e^{-4u} = K.
\end{equation}
Hence  all $f_{\lambda}: M \to \bb{S}^3_1$ are isometric surfaces. Summing up we have proved  the following

\begin{theorem}\label{S^1deformations}
Let $f: M \to \bb{S}^3_1$ be a conformal  immersion with constant mean curvature $H$,  gaussian curvature $K$ and Hopf complex differential $q =  \xi dz \otimes dz$.
Let $\widehat{f}:M \to \mathcal{Z}$ be its twistor lift. Then there is a one parameter family of harmonic maps $\phi_{\lambda} : M \to \mathcal{Z},\, \lambda \in \bb{S}^1$ satisfying  $\phi_{\{\lambda =1  \}} = \widehat{f}$, which are  
given by  $\phi_{\lambda} = F_{\lambda}.o$,  where the extended frame $F_{\lambda} : M \to SO_o(3,1)$ solves~\eqref{extendedframeeqs}.\\
The projection  $f_{\lambda} = \pi \circ \phi_{\lambda} : M \to \bb{S}^3_1$\,  is a $\bb{S}^1$-family of isometric conformally  immersed surfaces  satisfying  $f_{\{\lambda =1  \}} = f$ and $\widehat{f_{\lambda}} = \phi_{\lambda}, \,\, \forall \lambda \in \bb{S}^1$. Moreover, the induced metric  $f_{\lambda}^* \la.,.\ra$  does not depend on $\lambda$ and all  $f_{\lambda}$ have constant mean curvature $H$,   gaussian curvature $K$   and Hopf complex differential $q_{\lambda} = \lambda^{-2} \xi dz \otimes dz$.     
\end{theorem}

 \begin{remark}
  $\mathcal{Z}$ is a reductive homogeneous space which is  not  symmetric since  $[ \goth{p}, \goth{p} ] \not \subset \goth{k}$. Thus the harmonic map equation for arbitrary maps into $\mathcal{Z}$ cannot be encoded in a   loop of connections with zero curvature. Nevertheless since twistor lifts satisfy condition~\eqref{conditionX}  the harmonic map equation for twistor lifts admits a formulation  as the flatness condition (ZCC) of a family of connections  parameterized by unit complex numbers.  This reflects the complete integrability of the harmonic map equation for twistor lifts. 
 \end{remark}

\subsection{Holomorphic twistor lifts.} 
 Here  we  consider the behaviour of twistor lifts    in relation to  both invariant almost complex structures $J',J''$ on the horizontal distribution $\goth{h}\subset \mathcal{Z}$ introduced before.  
A  smooth map $ \phi :M \to \mathcal{Z}$ is said  {\it horizontal} if $d \phi (T_x M) \subset \goth{h}_{\phi(x)}$ for any $x \in M$. From the structure of matrices~\eqref{dzlift} and~\eqref{dzbarlift} we conclude that twistor lifts are horizontal maps.  Let $J $ be one of the almost complex structures $ J', J''$ considered before. A horizontal map $ \phi :M \to \mathcal{Z}$ is  $J$-holomorphic if it satisfies a Cauchy-Riemann type equation
\begin{equation}\label{JCRequation} 
J \circ d \phi = d \phi \circ J^M,
 \end{equation}
 where $J^M$ is the complex structure of $M$.    
Equivalenlty $\phi$ is $J$-holomorphic if and only if  $ d \phi (T^{(1,0)} M) \subset \goth{h}^{(1,0)}_{\phi}$, where 
$$\goth{h}^{(1,0)}_q = \{ X \in \goth{h}^{\C}_q : JX =iX \}.  $$

 Recall now that the isomorphism  $\beta : T \mathcal{Z} \to [\goth{p} ]$ constructed before satisfies \eqref{betaH}, i.e.      
  $$\beta_{g.o}(\goth{h}_{g.o}) = \{ g.o \} \times Ad(g) \mathcal{H},\, \forall g \in SO_o(3,1).$$
   Then  a horizontal map $\phi : M \to \mathcal{Z}$ is $J$-holomorphic if and only if  for every frame $F$ of $\phi$
  $$ 
       \phi^* \beta (\Pz) \in  Ad(F) \mathcal{H}^{(1,0)}.
  $$
   On the other hand for every frame $F$ of $\phi$   the following  identity holds 
  $$ \phi^* \beta (\Pz) = Ad(F) \alpha_{\goth{p}}(\Pz),$$ in which $\alpha_{\goth{p}}$ is the $\goth{p}$-component of the Maurer-Cartan one form $\alpha = F^{-1} dF$.    
  We  conclude that a horizontal map $\phi:M \to \mathcal{Z}$ is $J$-holomorphic if and only if  for every frame $F$ of $\phi$
  $$
          \alpha_{\goth{p}}(\Pz) \in \mathcal{H}^{(1,0)}. 
          $$

From the explicit form of $J'$ at $o \in \mathcal{Z}$  we see that the $i$-eigenspace $\mathcal{H}^{(1,0)}$ corresponding to $J'$  consists of matrices of the form 
$$
\begin{pmatrix}
0 &  a  & -ia & 0 \\
-a&  0  & 0  & b \\
ia& 0  &  0&  -ib \\
0 & b &   -ib &0  \\
\end{pmatrix}, a,b \in \C. $$ 

Now let $f: M \to \bb{S}^3_1$ be a conformal (hence spacelike) immersed surface and  $\widehat{f} :M \to \mathcal{Z}$ its twistor lift. 
 Thus $\widehat{f}$ is $J$-holomorphic if and only if  $\alpha_{\goth{p}}(\Pz) \in \mathcal{H}^{(1,0)}$ for every adapted frame $F$ of $f$. From~\eqref{dzlift} we have 
 $$   
 \alpha_{\goth{p}}(\Pz) = A_{\goth{p}}=
 \begin{pmatrix}
0   &   -\frac{e^u}{\sqrt{2}} & i\frac{e^u}{\sqrt{2}} & 0 \\
\frac{e^u}{\sqrt{2}}& 0 &  0 & \frac{ e^{-u} \xi + e^u H}{\sqrt{2}}\\
-i \frac{e^u}{\sqrt{2}} & 0 & 0 & i (  \frac{e^{-u} \xi-e^u H}{\sqrt{2}})\\
0 &   \frac{ e^{-u} \xi+ e^{u}H }{\sqrt{2}} &  i ( \frac{e^{-u}
\xi-e^u H}{\sqrt{2}})&
0 \\
\end{pmatrix},$$
  
which is in the $i$-eigenspace $\mathcal{H}^{(1,0)}$ of $J'$ if and only if    $i(e^{-u} \xi - e^{u} H ) = -i (e^{-u} \xi + e^{u} H )$,
 or  $2 e^{-u} \xi =0$, hence $\xi =0$.\\ 

Arguing in an analogous way we conclude that a horizontal map $\phi : M \to \mathcal{Z}$ is $J''$-holomorphic if and only if   for every frame $F$ of $\phi$  the complex matrix  $ \alpha_{\goth{p}}(\Pz)=A_{\goth{p}}$ is an eigenvector of $J''$ corresponding to the eigenvalue $i$.  Using the explicit form of $J''$ \eqref{J''}  we see that $X \in \mathcal{H}^{\C}$ satisfies $J'' X = iX$ if and only if  
$$
X = \begin{pmatrix}
0 & a & -ia & 0 \\
-a& 0 &  0  & b \\
ia& 0 &  0  & ib\\
0 & b & ib  & 0
\end{pmatrix}, a,b \in \C.
$$
 Turning to  the particular case of the twistor lift $\widehat{f}$ of an immersed surface $f$, we see that $\widehat{f}$ is $J''$-holomorphic if and only if     $J'' A_{\goth{p}} = i A_{\goth{p}}$ if and only if  $i(e^{-u} \xi - e^{u} H )= i(e^{-u} \xi + e^{u} H )$, or  $2i e^{u} H =0 $.  \\
 One can also  characterize twistor lifts which are {\it conformal} maps. For, let $f:M \to \bb{S}^3_1$ be a conformal immersion and take an  adapted (local) frame $F$ of  $f$. Thus $\widehat{f}^* \beta (\Pz)= Ad (F) A_{\goth{p}} $, where the complex matrix $A_{\goth{p}} $ is given by \eqref{dzlift}.
Since   $\beta$ preserves the metric,  we compute
\begin{equation*}
\begin{array}{cc}
\la  \widehat{f}_*(\Pz),\widehat{f}_*(\Pz) \ra^c=
\la \widehat{f}^* \beta (\Pz),\widehat{f}^* \beta (\Pz) \ra^c=
 \\ \la Ad(F)A_{\goth{p}}, Ad(F)A_{\goth{p}} \ra^c = 
 \la A_{\goth{p}}, A_{\goth{p}} \ra^c =  -\frac{1}{2} tr (A_{\goth{p}}^2) = -2 \xi H. \\
 \end{array} 
 \end{equation*} 
 Thus $\widehat{f}$ is conformal if and only if it is $J'$ or $J''$ holomorphic. 
 We have thus obtained the following.     

\begin{proposition}\label{holomorphictwistor}
Let $f:M \to \bb{S}^3_1$ be an immersed spacelike surface and let $\widehat{f}:M \to \mathcal{Z}$ its twistor lift. Then\\
i)  $\widehat{f}$ is $J'$-holomorphic if and only if $f$  is totally umbilic ($ \xi \equiv 0$).\\
ii) $\widehat{f}$ is $J''$-holomorphic if and only if $f$  has vanishing mean curvature ($H \equiv 0$).  \\
iii) $\widehat{f}$ is conformal if and only if $f$ satisfies $\xi.H =0$.  In particular if $\widehat{f}$ is conformal, then it is harmonic. 
\end{proposition}  
\begin{remark} As consequence of Theorem \ref{harmonictwsitorlift} and Codazzi's equation,  it follows that $J'$ and  $J''$-holomorphic twistor lifts  $\widehat{f}$ are harmonic maps. Moreover from the proof of Theorem~\ref{S^1deformations} it follows that the one parameter deformation introduced preserves $J'$ and $J''$-holomorphicity.   
\end{remark}

\section{ On the twistor energy}

  In order to gain some insight  of  the energy~\eqref{energyonomega} we compute the twistor energy and study its behaviour for (compact) genus zero and genus one spacelike surfaces. \\
Let $f: M \to \bb{S}^3_1$ be a conformal immersion, hence   the energy density $  \| d \widehat{f} \|^2$ of the twistor lift  is by definition  
 $$
 \| d \widehat{f} \|^2 (p)=  \la d \widehat{f}(e_1),d \widehat{f}(e_1) \ra + \la d \widehat{f}(e_2), d \widehat{f}(e_2)\ra,
 $$ 
 where $\{ e_1, e_2 \}$ is any orthonormal basis of $T_p M$. To compute $\| d \widehat{f} \|^2$ we use   the induced metric   $g=f^*\la\, , \, \ra$ which is conformal  and  locally given by  $g = 2 e^{2u} dz \otimes d \bar{z}$, where $u$ is the conformal factor. Thus    
 $$ \frac{1}{2} \| d \widehat{f} \|^2 =   e^{-2u} \la \widehat{f}_z, \widehat{f}_{\bar{z}} \ra^c.$$
 
 Now let $F$ be a   local adapted frame of $f$, thus  a frame of $\widehat{f}$ too. From the  identity $\widehat{f}^* \beta =  Ad(F)\alpha_{\goth{p}}$, in which  the complex matrix $A_{\goth{p}} = [F^{-1} F_z ]_{\goth{p}}$  is given by~\eqref{dzlift},and $B_{\goth{p}} = \overline{A_{\goth{p}}}$, we obtain 
 \begin{equation} \label{energydensityf}
 \begin{array}{cc}
 \\
 \frac{1}{2} \| d \widehat{f} \|^2 =  e^{-2u} \la \widehat{f}_z, \widehat{f}_{\bar{z}} \ra^c =
    e^{-2u} \la \widehat{f}^* \beta (\Pz),  \widehat{f}^*\beta(\Bz) \ra^c =  \\
    \\
   e^{-2u} \la A_{\goth{p}}, B_{\goth{p}} \, \ra^c  
   = - e^{-2u} \frac{1}{2} tr(A_{\goth{p}}.B_{\goth{p}} )=\\
   \\
    e^{-2u} ( e^{2u}(1-H^2) - e^{-2u} |\xi|^2 )=   \\
    \\
    1-H^2 - e^{-4u} |\xi|^2.
    \\
   \end{array}   
 \end{equation}
 If $\lambda_1, \lambda_2$ are the principal curvatures of the immersed surface, it is easily seen that   
 $$     e^{-4u} |\xi|^2 =  \frac{1}{4}(\lambda_1 - \lambda_2)^2.   $$    
 Therefore on a relatively compact domain $\Omega \subset M$ we obtain the following formula for the energy of $\widehat{f}$, 
 \begin{equation}\label{energyOmega}
     E_{\Omega}(\widehat{f})=  \int_{\Omega} [1-H^2-\frac{1}{4}(\lambda_1 - \lambda_2)^2]dA,
 \end{equation}
 where $dA$ is the area element of $(M,g)$.\\   
 \medskip

 On the other hand if $M$ is compact without boundary,  the Willmore energy of the conformal immersion  $f:M \to \bb{S}^3_1$ is given by
 \begin{equation}\label{willmoreenergy}
  W(f) = \int_M (K+H^2-1)dA = \frac{1}{4} \int_M (\lambda_1 - \lambda_2)^2 dA.
  \end{equation}
       
Combining~\eqref{energyOmega} and~\eqref{willmoreenergy} with Gauss equation~\eqref{gausseq2} and the Gauss-Bonnet formula we obtain
\begin{lemma} 
Let $f:M \to \bb{S}^3_1$ be a conformal immersion of a compact closed Riemann surface $M$. Then the  total energy of $\widehat{f}$ over $M$ and the Willmore energy $W(f)$ are related by the equality  
\begin{equation}\label{Willmore-E}
                   2W(f) = 2 \pi \mathcal{X}(M) - E(\widehat{f}),
\end{equation}
where $H$ and $K$ are  the mean curvature and the gaussian curvature respectively of the immersed surface, and $\mathcal{X}(M)$ is the Euler-Poincar\'e characteristic of $M$. 
\end{lemma}

 Since $W(f) \geq 0$ for every conformal immersion $f$, we deduce that 
 \begin{equation}\label{ineqE}
 E(\widehat{f}) \leq 2 \pi \mathcal{X}(M).
 \end{equation}
 
 If $f: M \to \bb{S}^3_1$ has constant mean curvature $H$ satisfying $H^2 <1$, Akutagawa~\cite{akutagawa} and independently   Ramanathan~\cite{ramanathan} proved that  $f(M)$ is a totally umbilic $2$-sphere with constant gaussian curvature $K=1-H^2 >0$. Thus $W(f)=0$ and so equality is attained in~\eqref{ineqE}, namely $E(\widehat{f}) = 4 \pi$. Thus there are no compact genus zero surfaces which are vacua of the twistor energy.\\
 
 On the other hand from~\eqref{Willmore-E} since $\mathcal{X}(T^2) =0$,  we obtain $E(\widehat{f}) \leq 0$ for every conformal immersion  $f: T^2 \to \bb{S}^3_1$  of the two torus $T^2 = \bb{S}^1 \times \bb{S}^1$. Thus  the mean curvature function $H$ of these immersions  satisfies $H^2 \geq 1$.\\ 
 One may wonder if there are  spacelike tori with zero twistor energy.    
Assume that  $f: T^2 \to \bb{S}^3_1$ is a  conformal immersion with  
  $E(\widehat{f}) = 0 $, then its  mean curvature should satisfy $H^2 =1$.  
 Lifting  $f$ to the universal covering $\bb{R}^2 \to T^2 = \bb{R}^2/\Gamma$ we obtain a double periodic (with respect to $\Gamma$) conformal immersion $\tilde{f} : \bb{R}^2 \to \bb{S}^3_1$ such that $\tilde{f}(\bb{R}^2)= f(T^2)$. The corresponding complex holomorphic function  $\tilde{\xi}$ is  entire and double-periodic, hence constant. If $f$ is umbilic free, then $ \tilde{\xi} \equiv 0$, and one can normalize (by a change of coordinate) so that  $\tilde{\xi} = 1$ on $\bb{R}^2$ and so $\xi =1$ on $T^2$.  Thus  the conformal immersion $f$ is free of umbilic points on the whole $T^2$,  so that  its principal curvatures never coincide  on $T^2$. In particular equation~\eqref{gausseq2} becomes  
  $$K = \frac{1}{4} (\lambda_1-\lambda_2)^2 >0$$    on the whole $T^2$.   Integrating this equation we  have
  $$   0 = \mathcal{X}(T^2)=  \int_{T^2} K dA =  W(f) >0.$$
  This shows that there is no genus one umbilic free spacelike surface in $\bb{S}^3_1$ with zero-energy twistor lift. On the other hand it is not difficult to rule out the existence of totally umbilic spacelike tori with $H^2=1$, i.e. vacua of the twistor energy. 
 Thus the search of genus $g > 1$ compact spacelike vacua of the twistor energy seems to be an interesting problem.   
 
 \begin{remark} Conformally immersed umbilic free  $f: T^2 \to \bb{S}^3_1$ with  constant mean curvature $H$ satisfying  $H^2 >1$ are determined by  double periodic solutions of  Gauss equation~\eqref{gauss-codazzi} which normalized is  the Sinh-Gordon equation: $$u_{\bar{z}z} + sinh(u) =0,$$ see~\cite{bobenko},~\cite{bobenko2},~\cite{palmer}. 
 \end{remark}

  \end{document}